\documentclass[12pt]{article}
\usepackage{amssymb}
\usepackage{mathrsfs}
\usepackage{pst-poly}     
\usepackage{cite}

\newtheorem{thm}{Theorem}[section]
\newtheorem{lem}[thm]{Lemma}
\newtheorem{Def}[thm]{Definition}
\newtheorem{cor}[thm]{Corollary}

\newenvironment{pf}[1][Proof]{\noindent\textbf{#1.} }{\hfill\rule{1mm}{2mm}}

\makeatletter \@addtoreset{equation}{section} \makeatother

\begin{document}
\title{\bf Generalized Measures of Edge Fault
Tolerance in $(n,k)$-star Graphs\thanks {The work was supported by
NNSF of China (No.11071233).}}
\author
{Xiang-Jun Li \quad Jun-Ming
Xu\footnote{Corresponding author: xujm@ustc.edu.cn (J.-M. Xu)}\\
{\small School of Mathematical Sciences, University of Science
and Technology of China,}\\
{\small  Wentsun Wu Key Laboratory of CAS, Hefei, 230026, China}  \\
 }
\date{}
 \maketitle

\begin{abstract}

This paper considers a kind of generalized measure $\lambda_s^{(h)}$
of fault tolerance in the $(n,k)$-star graph $S_{n,k}$ for $2
\leqslant k \leqslant n-1$ and $0\leqslant h \leqslant n-k$, and
determines $\lambda_s^{(h)}(S_{n,k})=\min\{(n-h-1)(h+1),
(n-k+1)(k-1)\}$, which implies that at least
$\min\{(n-k+1)(k-1),(n-h-1)(h+1)\}$ edges of $S_{n,k}$ have to
remove to get a disconnected graph that contains no vertices of
degree less than $h$. This result shows that the $(n,k)$-star graph
is robust when it is used to model the topological structure of a
large-scale parallel processing system.

\vskip6pt

\noindent{\bf Keywords:} Combinatorics, fault-tolerant analysis,
$(n,k)$-star graphs, edge-connectivity, $h$-super edge-connectivity

\end{abstract}

\section{Introduction}
It is well known that interconnection networks play an important
role in parallel computing/communication systems. An interconnection
network can be modeled by a graph $G=(V, E)$, where $V$ is the set
of processors and $E$ is the set of communication links in the
network. The connectivity $\lambda(G)$ of a graph $G$ is an
important measurement for fault-tolerance of the network, and the
larger $\lambda(G)$ is, the more reliable the network is.

A subset of vertices $B$ of a connected graph $G$ is called a {\it
edge-cut} if $G-B$ is disconnected. The {\it edge connectivity}
$\lambda(G)$ of $G$ is defined as the minimum cardinality over all
edge-cuts of $G$. Because $\lambda$ has many shortcomings, one
proposes the concept of the $h$-super edge connectivity of $G$, which can
measure fault tolerance of an interconnection network more
accurately than the classical connectivity $\lambda$.

A subset of vertices $B$ of a connected graph $G$ is called an {\it
$h$-super edge-cut}, or {\it $h$-edge-cut} for short, if $G-B$ is
disconnected and has the minimum degree at least $h$. The {\it
$h$-super edge-connectivity} of $G$, denoted by
$\lambda^{(h)}_s(G)$, is defined as the minimum cardinality over all
$h$-edge-cuts of $G$. It is clear that, if $\lambda_s^{(h)}(G)$
exists, then
 $$
 \lambda(G)=\lambda_s^{(0)}(G)\leqslant \lambda_s^{(1)}(G)\leqslant \lambda_s^{(2)}(G)\leqslant
 \cdots \leqslant  \lambda_s^{(h-1)}(G)\leqslant \lambda_s^{(h)}(G).
 $$

For any graph $G$ and integer $h$, determining $\lambda_s^{(h)}(G)$
is quite difficult. In fact, the existence of $\lambda_s^{(h)}(G)$
is an open problem so far when $h\geqslant 1$. Some results have
been obtained on $\lambda_s^{(h)}$ for particular classes of graphs
and small $h$'s (see Section 16.7 in~\cite{x07}).

This paper is concerned about $\lambda_s^{(h)}$ for the $(n,k)$-star
graph $S_{n,k}$. In $h$-super connectivity, several authors
 have done some work. For $k=n-1$, $S_{n,n-1}$ is isomorphic to a star
graph $S_n$. Akers and Krishnamurthy~\cite{ak89} determined
$\lambda(S_n)=n-1$ for $n\geqslant 2$ and
$\lambda_s^{(1)}(S_n)=2n-4$ for $n\geqslant 3$. In this paper, we
show the following result.

\begin{thm}\label{thm1.1}
If $2 \leqslant k \le n-1$ and $0\leqslant h \leqslant n-k$, then
 $$
 \lambda_s^{(h)}(S_{n,k})= \left\{\begin{array}{ll}
 (n-h-1)(h+1)&{}\  for\ h\leqslant k-2  \ and \ h\leqslant \frac{n}{2}-1,\\
 (n-k+1)(k-1)&{}\  otherwise.\\
    \end{array}
\right.$$
\end{thm}

This result implies that at least
$\min\{(n-k+1)(k-1),(n-h-1)(h+1)\}$ edges of $S_{n,k}$ have to
remove to get a disconnected graph that contains no vertices of
degree less than $h$. The proof of this result is in Section 3. In
Section 2, we recall the structure of $S_{n,k}$ and some lemmas used
in our proofs.

\section{Definitions and lemmas}

For given integer $n$ and $k$ with $1\leqslant k\leqslant n-1$, let
$I_n=\{1,2,\ldots,n\}$ and $P(n,k)=\{ p_{1}p_{2}\ldots p_{k}:\
p_{i}\in I_n, p_{i}\neq p_{j}, 1\leqslant i\neq j\leqslant k\}$, the
set of $k$-permutations on $I_n$. Clearly, $|P(n,k)|=n\,!/(n-k)\,!$.

\begin{Def}\label{def2.1}
The $(n,k)$-star graph $S_{n,k}$ is a graph with vertex-set
$P(n,k)$. The adjacency is defined as follows: a vertex
$p=p_{1}p_{2}\ldots p_{i}\ldots p_{k}$ is adjacent to a vertex

(a)\ $p_{i}p_{2}\cdots p_{i-1}p_{1}p_{i+1}\cdots p_{k}$, where
$2\leqslant i\leqslant k$ (swap $p_{1}$ with $p_{i}$).

(b)\ $\alpha p_{2}p_{3}\cdots p_{k}$, where $\alpha\in I_n\setminus
\{p_{i}:\ 1\leqslant i\leqslant k\}$ (replace $p_{1}$ by $\alpha$).
\end{Def}

The vertices of type $(a)$ are referred to as {\it swap-neighbors}
of $p$ and the edges between them are referred to as {\it swap-edge}
or {\it $i$-edges}. The vertices of type $(b)$ are referred to as
{\it unswap-neighbors} of $p$ and the edges between them are
referred to as {\it unswap-edges}. Clearly, every vertex in
$S_{n,k}$ has $k-1$ swap-neighbors and $n-k$ unswap-neighbors.
Usually, if $x=p_1p_2\dots p_k$ is a vertex in $S_{n,k}$, we call
$p_i$ the {\it $i$-th bit} for each $i\in I_k$.

The $(n,k)$-star graph $S_{n,k}$ is proposed by Chiang and
Chen~\cite{cc95}. Some nice properties of $S_{n,k}$ are compiled by
Cheng and Lipman (see Theorem 1 in~\cite{cl02a}).

\begin{lem}\label{lem2.2}
$S_{n,k}$ is $(n-1)$-regular $(n-1)$-connected.
\end{lem}

\begin{lem}\label{lem2.3}
For any $\alpha=p_1p_2\cdots p_{k-1}\in P(n,k-1)$ $(k \geqslant 2)$,
let $V_\alpha=\{p\alpha:\ p\in I_n\setminus \{p_i:\ i\in
I_{k-1}\}\}$. Then the subgraph of $S_{n,k}$ induced by $V_\alpha$
is a complete graph of order $n-k+1$, denoted by $K^\alpha_{n-k+1}$.
\end{lem}

Let $S^{t:i}_{n-1,k-1}$ denote a subgraph of $S_{n,k}$ induced by
vertices with the $t$-th bit $i$ for $2\leqslant t\leqslant k$. The
following lemma is a slight modification of the result of Chiang and
Chen~\cite{cc95}.

\begin{lem}\label{lem2.4}
For a fixed integer $t$ with $2\leqslant t\leqslant k$, $S_{n,k}$
can be decomposed into $n$ subgraphs $S^{t:i}_{n-1,k-1}$, which is
isomorphic to $S_{n-1,k-1}$, for each $i\in I_n$. Moreover, there
are $\frac{(n - 2)!}{(n - k)!}$ independent swap-edges between
$S^{t:i}_{n-1,k-1}$ and $S^{t:j}_{n-1,k-1}$ for any $i,j\in I_n$
with $i\ne j$.
\end{lem}

Since $S_{n,1}\cong K_n$, we only consider the case of $k\geqslant
2$ in the following discussion.

\begin{lem}\label{lem2.5}
If $2 \leqslant k \le n-1$ and $0\leqslant h \leqslant n-k$, then
 $$
 \lambda_s^{(h)}(S_{n,k})\leqslant \left\{\begin{array}{ll}
 (n-h-1)(h+1)&\  {\rm for}\  \ h\leqslant \frac{n}{2}-1,\\
 (n-k+1)(k-1)&\  {\rm otherwise}.\\
    \end{array}
\right.$$
\end{lem}

\begin{pf}
By our hypothesis of $h \leqslant n-k$, for any $\alpha\in
P(n,k-1)$, we can choose a subset $X\subseteq V(K^\alpha_{n-k+1})$
such that $|X|=h+1$. Then the subgraph of $K^\alpha_{n-k+1}$ induced
by $X$ is a complete graph $K_{h+1}$. Let $B$ be the set of incident
edges with and not within $X$. Since $S_{n,k}$ is $(n-1)$-regular and
$K_{h+1}$ is $h$-regular, we have that
 $$
 |B|=(n-h-1)(h+1).
 $$

 Clearly, $B$ is an edge-cut of $S_{n,k}$.  Let $x$ be any vertex in$S_{n,k}-X$,
 and $d(x)$ denote the
number of edges incident with $x$ in $S_{n,k}-X$. In order to prove that $B$ is an $h$-edge-cut, we
only need to show $d(x) \geqslant h$.
Note that $X$ is contained in $S^i_{n-1,k-1}$
and edges between $S^i_{n-1,k-1}$  and $S^j_{n-1,k-1}$  are
independent for any $i, j \in I_n$ with $i\not= j$ by Lemma 2.4.
If $x$ is in $S^i_{n-1,k-1}-K_{n-k+1}$ or is in $S^j_{n-1,k-1}$
with $i\not= j$, then $d(x) \geqslant n - 2 \geqslant n - k \geqslant h$. For $x \in V (K_{n-k+1}-X)$, if exists, then
$d(x) = n-1-|X| = n- h-2 \geqslant h$ for $h \leqslant \frac{n}{2}-1$.
Therefore, $B$ is an $h$-edge-cut of $S_{n,k}$, and so
 $$
 \lambda_s^{(h)}(S_{n,k})\leqslant |B|=(n-h-1)(h+1) \ \ {\rm for } \ h\leqslant \frac{n}{2}-1.
 $$

If $h\geqslant\frac{n}{2}$, we choose $X=V(K^\alpha_{n-k+1})$.
Then $|B|=(n-k+1)(k-1)$. For any $x$ in $S^i_{n-1,k-1}-X$ or
$S^j_{n-1,k-1}$ with $i\not= j$, we have $d(x) \geqslant n -
2\geqslant n - k \geqslant h$. Thus, $B$ is an $h$-edge-cut of
$S_{n,k}$, and so
 $$
 \lambda_s^{(h)}(S_{n,k})\leqslant |B|=(n-k+1)(k-1) \ \ {\rm for } \ h\geqslant\frac{n}{2}.
 $$
The lemma follows.
\end{pf}

\begin{cor}\label{cor2.6}
$\lambda_s^{(h)}(S_{n,2})=n-1$ for $0 \leqslant h\leqslant n-2$.
\end{cor}

\begin{pf}
On the one hand, $\lambda _s^{(h)} (S_{n,2})\leqslant n-1$ by
Lemma~\ref{lem2.5} when $k=2$. On the other hand, $\lambda _s^{(h)}
(S_{n,2})\ge \lambda(S_{n,2})=n-1$ by Lemma~\ref{lem2.2}.
\end{pf}

\vskip6pt The following lemma shows the relations between
$(n-h-1)(h+1)$ and $(n-k+1)(k-1)$.

\begin{lem}\label{lem2.7}
For $2\leqslant k\leqslant n-1,0 \leqslant h \leqslant n-k$, let
 \begin{equation}\label{e3.0a}
 \psi(h,k)=\min\{(n-h-1)(h+1),(n-k+1)(k-1)\}.
 \end{equation}
If $h\leqslant \frac n2-1$, then
  $$
 \psi(h,k)=\left \{\begin{array}{rl}
 (n-h-1)(h+1)& \ {\rm if}  \ 0\leqslant h\leqslant k-2;\\
 (n-k+1)(k-1)& \ {\rm if}   \ h\geqslant k-1.
 \end{array}\right.
 $$
\end{lem}

\begin{pf}
Let $f(x)=(n-x)\,x$, then $\psi(h,k)=\min\{f(h+1),f(k-1)\}$. It can
be easily checked that $f(x)$ is a convex function on the interval
$[0,n]$, the maximum value is reached at $x=\frac n2$. Thus, $f(x)$
is an increasing function on the interval $[0,\frac n2]$.

If $ 0\leqslant h\leqslant k-2$, then $h+1\leqslant k-1$. Since
$h\leqslant n-k$, $h+1 \leqslant n-k+1$. Clearly, $\max\{k-1,
n-k+1\}\leqslant\frac n2$. Thus, when $h\leqslant \frac n2-1$,
$f(h+1)\leqslant f(k-1)=f(n-k+1)$, and so
$\psi(h,k)=f(h+1)=(n-h-1)(h+1)$.

If $h\geqslant k-1$, then $k-1<h+1\leqslant \frac {n}2$. Thus,
$f(k-1)<f(h+1)$, and so $\psi(n,k)=f(k-1)=(n-k+1)(k-1)$.

The lemma follows.
\end{pf}

\vskip6pt

To state and prove our main results, we need some notations. Let $B$
be a minimum $h$-edge-cut of $S_{n,k}$. Clearly, $S_{n,k}-B$ has
exactly two connected components. Let $X$ and $Y$ be two vertex-set
of two connected components of $S_{n,k}-B$. For a fixed $t\in
I_k\setminus\{1\}$ and any $i\in I_n$, let
 \begin{equation}\label{e3.2}
 \begin{array}{rl}
 &X_i=X\cap V(S^{t:i}_{n-1,k-1}),\\
 &Y_i=Y\cap V(S^{t:i}_{n-1,k-1}),\\
 &B_i=B\cap E(S^{t:i}_{n-1,k-1})\ {\rm and}\\
 &B_{ij}=B\cap E(S^{t:i}_{n-1,k-1},S^{t:j}_{n-1,k-1}),
 \end{array}
\end{equation}
and let
\begin{equation}\label{e3.3}
\begin{array}{l}
 J=\{i\in I_n:\ X_i\ne\emptyset\},\\
 J'=\{i\in  J:\ Y_i\not=\emptyset\}\ \ {\rm and} \\
 T=\{i\in I_n:\  Y_i\ne\emptyset\}.
 \end{array}
 \end{equation}

\begin{lem}\label{lem2.8}
Let $B$ be a minimum $h$-edge-cut of $S_{n,k}$ and $X$ be the
vertex-set of a connected component of $S_{n,k}-B$. If $3\leqslant
k\leqslant n-1$ and $1\leqslant h\leqslant n-k$ then, for any $t\in
I_k\setminus\{1\}$,

{\rm (a)}\ $B_i$ is an $(h-1)$-edge-cut of $S^{t:i}_{n-1,k-1}$ for
any $i\in J'$,

{\rm (b)}\ $\lambda_s^{(h)}(S_{n,k})\geqslant |J'|\
\lambda_s^{(h-1)}(S_{n-1,k-1})$,

\end{lem}

\begin{pf}
(a)\ By the definition of $J'$, $B_i$ is an edge-cut of
$S^{t:i}_{n-1,k-1}$ for any $i\in J'$. For any vertex $x$ in
$S^{t:i}_{n-1,k-1}-B_i$, since $x$ has degree at least $h$ in
$S_{n,k}-S$ and has exactly one neighbor outsider
$S^{t:i}_{n-1,k-1}$, $x$ has degree at least $h-1$ in
$S^{t:i}_{n,k}-B_i$. This fact shows that $B_i$ is an
$(h-1)$-edge-cut of $S^{t:i}_{n-1,k-1}$ for any $i\in J'$.

(b)\ By the assertion (a), we have
$|B_i|\geqslant \lambda_s^{(h-1)}(S_{n-1,k-1})$, and so
 $$
 \lambda_s^{(h)}(S_{n,k})=|B|\geqslant \sum_{i\in J'} |B_i|\geqslant
 |J'|\lambda_s^{(h-1)}(S_{n-1,k-1}).
  $$
the lemma follows.
\end{pf}

\section{Proof of Theorem~\ref{thm1.1}}

\begin{pf}
By Lemma~\ref{lem2.5} and  Lemma~\ref{lem2.7}, we only need to prove
that, for $2\leqslant k\leqslant n-1$ and $0\leqslant h\leqslant
n-k$,
 \begin{equation}\label{e3.4}
 \lambda_s^{(h)}(S_{n,k})\geqslant \left\{\begin{array}{ll}
 (n-h-1)(h+1)&{}\ {\rm for}\ h\leqslant k-2 \ {\rm and} \ h\leqslant \frac{n}{2}-1,\\
 (n-k+1)(k-1)&{}\  {\rm otherwise}.\\
    \end{array}
\right.
 \end{equation}

Let $\omega(h,k)=\max\{(n-h-1)(h+1),(n-k+1)(k-1)\}$.

We proceed by induction on $k\geqslant 2$ and $h\geqslant 0$. The
inequality (\ref{e3.4}) is true for $k=2$ and any $h$ with $0
\leqslant h\leqslant n-2$ by Corollary~\ref{cor2.6}. The inequality
(\ref{e3.4}) is also true for $h=0$ and any $k$ with $2\leqslant
k\leqslant n-1$ since $\lambda
_s^{(0)}(S_{n,k})=\lambda(S_{n,k})=n-1$. Assume the induction
hypothesis for $k-1$ with $k\geqslant 3$ and for $h-1$ with
$h\geqslant 1$, that is,
 \begin{equation}\label{e3.5}
 \lambda_s^{(h-1)}(S_{n-1,k-1})\geqslant \left\{\begin{array}{ll}
 (n-h)h&{}\ {\rm for}\ h\leqslant k-3 \ {\rm and} \ h\leqslant \frac{n-1}{2},\\
 (n-k+2)(k-2)&{}\ {\rm otherwise}.\\
    \end{array}\right.
 \end{equation}

Let $B$ be a minimum $h$-edge-cut of $S_{n,k}$ and $X$ be the
vertex-set of a minimum connected component of $S_{n,k}-B$. By
Lemma~\ref{lem2.5}, we have
 \begin{equation}\label{e3.5a}
 |B|\leqslant \omega(h,k).
\end{equation}

Use notations defined in (\ref{e3.2}) and (\ref{e3.3}). Choose $t\in
I_k\setminus\{1\}$ such that $|J|$ is as large as possible. For each
$i\in I_n$, we write $S^i_{n-1,k-1}$ for $S^{t:i}_{n-1,k-1}$ for
short.

We first show $|J|=1$. Suppose to the contrary $|J|\geqslant 2$. We
will deduce contradictions by considering three cases depending on
$|J'|=0$, $|J'|=1$ or $|J'|\geqslant2$.

\vskip6pt

{\bf Case 1.}\  $|J'|=0$,

In this case, $X_i\ne\emptyset$ and $Y_i=\emptyset$ for each $i\in
J$, that is, $J\cap T=\emptyset$. By $|J|\geqslant2 $ and the
minimality of $X$, $|T|\geqslant 2$. Assume $\{i_{1},i_2\}\subseteq
J$ and $\{i_3,i_4\}\in T$. By Lemma~\ref{lem2.4}, there are
$\frac{(n - 2)!}{(n - k)!}$ independent swap-edges between
$S^{i_1}_{n-1,k-1}$ (resp. $S^{i_2}_{n-1,k-1}$) and
$S^{i_3}_{n-1,k-1}$(resp. $S^{i_4}_{n-1,k-1}$), all of which are
contained in $B$. Since $J\cap T=\emptyset$, we have that
$$
\begin{array}{c}
 |B|\geqslant 4\ \frac{(n-2)!}{(n-k)!}.
  \end{array}
$$
For $k=3$,
 $$
 \begin{array}{c}
|B|\geqslant 4\ \frac{(n-2)!}{(n-k)!}\geqslant 4(n-2)>2(n-2)
  \end{array}
 $$
Combining Lemma~\ref{lem2.5} with Lemma~\ref{lem2.7} yields $|B|
\leqslant \lambda_s^{(h)}(S_{n,3}) \leqslant 2(n-2)$, a
contradiction. For $k\geqslant 4$, it is easy to check that
 $$
 \begin{array}{rl}
 |B|&\geqslant 4\ \frac{(n-2)!}{(n-k)!}\geqslant 4(n-2)(n-3)=(2n-4)(2n-6)\\
  &>\max\{(n-h-1)(h+1),(n-k+1)(k-1)\}\\
  &= \omega(h,k),
  \end{array}
 $$
which contradicts the inequality (\ref{e3.5a}).

\vskip6pt

{\bf Case 2.}\  $|J'|=1$,

Without loss of generality, assume $J'=\{1\}$. By Lemma~\ref{lem2.8}
(a), $B_1$ is an $(h-1)$-edge-cut of $S^1_{n-1,k-1}$.

By $|J|\geqslant2$, there exists an $i\in J-J'$ such that
$X_i=V(S^i_{n-1,k-1})$. By the minimality of $X$, there exists some
$j\in T-J'$ such that $Y_j=V(S^j_{n-1,k-1})$.  By
Lemma~\ref{lem2.4}, there are $\frac{(n - 2)!}{(n - k)!}$
independent swap-edges between $S^{i}_{n-1,k-1}$ and
$S^{j}_{n-1,k-1}$, thus $|B_{ij}|=\frac{(n - 2)!}{(n - k)!}\geqslant
n-2$. By (\ref{e3.5}), we consider the following two cases.

If $\lambda_s^{(h-1)}(S_{n-1,k-1})\geqslant (n-h)h$, then
 $$\begin{array}{rl}|B|&\geqslant |B_1|+ |B_{ij}|\\
 &\geqslant (n-h)h+(n-2)\\&>(n-h-1)h+(n-h-1)\\
 &=(n-h-1)(h+1),\end{array}$$
If $\lambda_s^{(h-1)}(S_{n-1,k-1})\geqslant (n-k+2)(k-2)$, then
 $$\begin{array}{rl}|B|&\geqslant |B_1|+ |B_{ij}|\\
  &\geqslant (n-k+2)(k-2)+(n-2)\\&>(n-k+2)(k-2)+(n-k+2)\\&>(n-k+1)(k-1).\end{array}$$
Therefore, we have $|B|>\omega(h,k)$, which contradicts the
inequality (\ref{e3.5a}).

{\bf Case 3.}\  $|J'|\geqslant 2$.

By Lemma~\ref{lem2.8} (b) and (\ref{e3.5}), we consider the
following two cases.

If $\lambda_s^{(h-1)}(S_{n-1,k-1})\geqslant (n-h)h$, then
 $$
 \begin{array}{rl}
|B|& \geqslant |J'|\lambda_s^{(h-1)}(S_{n-1,k-1})\\
 &\geqslant 2(n-h)h\geqslant(n-h)h+(n-h)\\
 &>(n-h-1)(h+1),
 \end{array}
 $$
If $\lambda_s^{(h-1)}(S_{n-1,k-1})\geqslant (n-k+2)(k-2)$, then
 $$
 \begin{array}{rl}
|B| &\geqslant |J'|\lambda_s^{(h-1)}(S_{n-1,k-1})\\
 &\geqslant 2(n-k+2)(k-2)\\
 &\geqslant(n-k+2)(k-2)+(n-k+2)\\&>(n-k+1)(k-1).\end{array}$$
Therefore, we have $|B|>\omega(h,k)$, which contradicts the
inequality (\ref{e3.5a}).

Thus, we have $|J|=1$. By the choice of $t$, the $i$-th bits of all
vertices in $X$ are same for each $i=2,3,\ldots,k$, and so $X$ is a
complete graph. Thus, we have that
     $$
\begin{array}{rl}
\lambda_s^{(h)}(S_{n,k})& =|B|=(n-|X|)|X|
\end{array}
    $$

Since $h+1\leqslant |X|\leqslant n-k+1$ and $f(x)=(n-x)\,x$ is a
convex function on the interval $[0,n]$, we have that
 $$
\begin{array}{rl}
\lambda_s^{(h)}(S_{n,k})& =|B|=(n-|X|)|X|\geqslant \psi(h,k),
\end{array}
   $$
where $\psi(h,k)$ is defined in (\ref{e3.0a}).

If $ h\leqslant \frac{n}{2}-1$, using Lemma~\ref{lem2.7},
 we have
 \begin{equation}\label{e3.7}
 \lambda_s^{(h)}(S_{n,k})\geqslant \psi(h,k)=\left \{\begin{array}{rl}
 (n-h-1)(h+1)& \ {\rm if}  \ 0\leqslant h\leqslant k-2;\\
 (n-k+1)(k-1)& \ {\rm if}   \ h\geqslant k-1.
 \end{array}\right.
 \end{equation}

If  $h\geqslant\frac{n}{2}$, we have $X=V(K_{n-k+1})$. Otherwise,
there exists some $x\in V(K_{n-k+1}-X)$ such that
 $$
 h\leqslant d(x)=n-1-|X|\leqslant n-h-2,
 $$
which implies $h\leqslant \frac{n}{2}-1$, a contradiction.
Therefore, we have $|X|=n-k+1$, and
\begin{equation}\label{e3.8}
\begin{array}{rl}
\lambda_s^{(h)}(S_{n,k})=|B|=(n-|X|)|X|= (n-k+1)(k-1)\ {\rm for}\
h\geqslant\frac{n}{2}.
\end{array}
\end{equation}

Combining (\ref{e3.7}) with (\ref{e3.8}) yields (\ref{e3.4}). By the
induction principle, the theorem follows.
\end{pf}

\vskip6pt

As we have known, when $k=n-1$, $S_{n,n-1}$ is isomorphic to the
star graph $S_n$. Akers and Krishnamurthy~\cite{ak89} determined
$\lambda(S_n)$ and $\lambda_s^{(1)}(S_n)$, which can be obtained
from our result by setting $k=n-1$ and $h=0,1$, respectively.

\begin{cor}\textnormal{(Akers and Krishnamurthy~\cite{ak89})}
$\lambda(S_n)=n-1$ for $n\geqslant 2$ and
$\lambda_s^{(1)}(S_n)=2n-4$ for $n\geqslant 3$.
\end{cor}

\end{document}